\documentclass[english]{amsart}
\usepackage{amssymb}
\usepackage{amsmath}
\usepackage{amscd}
\usepackage{amsthm}
\usepackage{amsfonts}
\usepackage{graphicx}
\usepackage[all]{xy}

\theoremstyle{plain}
\newtheorem{stam}{STAM}[section]
\newtheorem{thm}[stam]{Theorem}

\newtheorem{definition}[stam]{Definition}
\numberwithin{equation}{section} 

\theoremstyle{plain}

\theoremstyle{plain}

\theoremstyle{remark}

\DeclareMathOperator{\aff}{aff}

\DeclareMathOperator{\conv}{conv}
\DeclareMathOperator{\rank}{rank}

\def\bbr{\mathbb{R}}

\begin{document}

\title{Strong general position}
\large

\author{Micha A. Perles, Moriah Sigron}

\begin{abstract}
We say that a finite set $S$  of points in $\bbr^d$ is in \textbf{strong general position} if  for any collection $\{ F_1, \ldots, F_r\}$ of $r$ pairwise
disjoint subsets of $S \quad (1 \le r \le |S|)$ we have:
 $d-\dim\bigcap^r_{\nu=1} \aff F_\nu = \min(d+1,\sum^r_{\nu=1}(d-\dim\aff F_\nu))$. In this paper we reduce the set of conditions that one has to check in order to determine if $S$  is  in \textbf{strong general position}.

\end{abstract}

\maketitle
\newsavebox{\rb}
\savebox{\rb}(10,40)[bl]{
  \put(0,2.5){\line(0,1){40}}
  \put(0,0){\circle{5}}
  \put(-2.5,10){\line(1,0){5}}
  \put(3,7){$b$}
  \put(-2.5,35){\line(1,0){5}}
  \put(3,32){$r$}
}
\newsavebox{\br}
\savebox{\br}(10,40)[bl]{
  \put(0,2.5){\line(0,1){40}}
  \put(0,0){\circle{5}}
  \put(-2.5,10){\line(1,0){5}}
  \put(3,7){$r$}
  \put(-2.5,35){\line(1,0){5}}
  \put(3,32){$b$}
}
\newsavebox{\ry}
\savebox{\ry}(10,40)[bl]{
  \put(0,2.5){\line(0,1){40}}
  \put(0,0){\circle{5}}
  \put(-2.5,10){\line(1,0){5}}
  \put(3,7){$y$}
  \put(-2.5,35){\line(1,0){5}}
  \put(3,32){$r$}
}
\newsavebox{\by}
\savebox{\by}(10,40)[bl]{
  \put(0,2.5){\line(0,1){40}}
  \put(0,0){\circle{5}}
  \put(-2.5,10){\line(1,0){5}}
  \put(3,7){$y$}
  \put(-2.5,35){\line(1,0){5}}
  \put(3,32){$b$}
}
\newsavebox{\yb}
\savebox{\yb}(10,40)[bl]{
  \put(0,2.5){\line(0,1){40}}
  \put(0,0){\circle{5}}
  \put(-2.5,10){\line(1,0){5}}
  \put(3,7){$b$}
  \put(-2.5,35){\line(1,0){5}}
  \put(3,32){$y$}
}
\newsavebox{\yr}
\savebox{\yr}(10,40)[bl]{
  \put(0,2.5){\line(0,1){40}}
  \put(0,0){\circle{5}}
  \put(-2.5,10){\line(1,0){5}}
  \put(3,7){$r$}
  \put(-2.5,35){\line(1,0){5}}
  \put(3,32){$y$}
}
\newsavebox{\simplebar}
\savebox{\simplebar}(10,40)[bl]{
  \put(0,2.5){\line(0,1){40}}
  \put(0,0){\circle{5}}
}
\newsavebox{\rightmost}
\savebox{\rightmost}(10,40)[bl]{
  \put(0,2.5){\line(0,1){40}}
  \put(0,0){\circle{5}}
  \put(-2.5,35){\line(1,0){5}}
  \put(3,32){$p_{5,1}$}
}

\setcounter{tocdepth}{1}

\section{Introduction}

A set $S$ of  points in $\bbr^d$ is said to be in \textbf{ general position} if every set of  $d+1$ or fewer points of $S$ is affinely independent, or, in other words, if any $k$-flat in $\bbr^d$ ($1 \leq k \leq d-1$) contains at most $k+1$ points of $S$. This condition can be simplified:

$S \subset \bbr^d $ is in general position if either $|S| \leq d+1$ and $S$ is affinely independent, or  $|S| \geq d+1$ and every $d+1$ points of $S$ are affinely independent.

General position is a rather weak property. E.g., the vertices  of a regular $2m$-gon $P  \subset  \bbr^2$ are in general position, even though the main diagonals of $P$ cross in a common point, and, moreover, opposite edges span parallel lines.

In this paper we consider a stronger property, called  \textbf{ strong general position} (SGP).  A finite set  $S\subset R^d$ is said to be in SGP if,  for any collection $\{ F_1, \ldots, F_r\}$ of $r$ pairwise
disjoint subsets of $S \quad (1 \le r \le |S|)$, the affine hulls $\aff F_1,\ldots,\aff F_r$ intersect as if they were flats chosen at random. The formal condition is :
\begin{equation}
\label{equ:Det_SGP1}
\dim\bigcap^r_{\nu=1} \aff F_\nu = \max(-1,
d-\sum^r_{\nu=1}(d-\dim\aff F_\nu)).
\end{equation}

As we shall see, this condition implies (ordinary) general position.

If we want to check whether a given large set $S \subset \bbr^d $ is in SGP, we are faced with a huge number of conditions of the form $(\ref{equ:Det_SGP1})$. The purpose of this paper is to reduce this number, i.e., to find a much smaller, essentially minimal set of conditions that will ensure that a given finite set $S \subset  \bbr^d $ is in SGP. This reduction  could be of use to anyone who wishes to work seriously with the notion of SGP. We use it in  \cite{PS} to show that points chosen on the moment curve $M^d$ in $ \bbr^d$ ($M^d = \{(t,t^2,\ldots,t^d) : t \in \bbr\}$)  are "usually" in SGP.

The notion of strong general position has been used (under the name "strong independence") by Reay in \cite{R} and by Doignon and Valette in  \cite{DV}.

Strong general position plays an important role in connection with Tverberg's  Theorem:

\begin{thm}
\label{thm:Tverberg}
(H. Tverberg, 1966) Let $a_{1},\ldots,a_{n}$ be points in $\mathbb{R}^{d}$.
If  $ n>(d+1)(r-1)$, then the set $ N=\{1,\ldots,n\}$ of indices
can be partitioned into r disjoint parts $N_{1},\ldots,N_{r}$ in such
a way that the r convex hulls $\conv \{ a_{i}:i\in N_{j}\}    (j=1,\ldots,r)$
have a point in common.
\end{thm}
(This formulation covers also the case where the points $a_{1},\ldots,a_{n}$
are not all distinct.)  The original proof (see \cite{T66}) was quite difficult.
In 1981 Tverberg published another proof, much simpler than the original one (see \cite{T81}).  Sarkaria \cite{Sa} gave a quite accessible proof, with some algebraic flavor. It seems that the simplest proof so far is due to Roudneff \cite{Ro}. See \cite{M} \S8.3 for further information.

The numbers $T(d,r)=(d+1)(r-1)+1$ are known as Tverberg numbers.  The condition $n\geq T(d,r)$ in Tverberg's theorem is extremely tight.  If $n<T(d,r)$, and the points $a_{1},\ldots,a_{n}$ are in SGP, then for any $r$-partition $N_{1},\ldots,N_{r}$ of the set $N=\{1,\ldots,n\}$, even the intersection of the \textbf{affine} hulls $\aff\{a_i:i \in N_j \} \,(j=1,\ldots,r)$ is empty. (See details in the next section.)

Our reduction (see Theorem $\ref{thm:SGP_result}$ below) will show that if $S \subset \bbr^d$ is a finite set in general position, then $S$ is in SGP iff, for any collection $\{F_1,\ldots,F_r \}$ of pairwise disjoint non-empty subsets of $S$  (with $2 \leq r \leq d+1,  |F_\nu| \leq d$ for all $i)$ of total size $m=\sum_{\nu=1}^r |F_\nu|$, the intersection $\cap_{i=1}^r \aff F_\nu$ is a single point if $m=T(d,r)$, or empty if   $m<T(d,r)$.

As we shall see in Section $\ref{sec:PointAreUsuallyInSGP}$, finite subsets of $\bbr^d$ are "usually" in SGP, in the following strong sense: Given $d$ and $n$, there exists a polynomial $P (=P_{d,n})$, not identically zero, in $nd$ scalar variables:  $P(\vec{x_1},\ldots,\vec{x_n})=P(x_{11},\ldots,x_{1d},\ldots,x_{n1},\ldots,x_{nd})$, such that any $n$ points  $\vec{a_{1}},\ldots,\vec{a_{n}} \in \bbr^d$ are (distinct and) in SGP unless $P( \vec{a_{1}},\ldots,\vec{a_{n}})=0$.

There are notions of independence that are even stronger than SGP. In fact, the first proof of Tverberg's Theorem in  \cite{T66} runs under the assumption that the points  $a_{1},\ldots,a_{n} \in \bbr^d$ are\textbf{ algebraically independent}, i.e., that the $nd$ coordinates $a_{i j}  (1 \leq i \leq n,  1 \leq j \leq d)$ are  algebraically independent over the field of rational numbers. A limiting argument then establishes  Tverberg's Theorem for all  $a_{1},\ldots,a_{n} \in \bbr^d$.
\section{Strong General Position}
\label{sec:SGP}

A $(d-k)$-dimensional flat in $\bbr^d \;(0\le k \le d)$ is the set
of solutions of a system of $k$ linearly independent linear
equations (not necessarily homogeneous) in $d$ variables.  It
follows that if $A_1,\ldots, A_r$ are flats in $\bbr^d$, and $\dim
A_\nu = d -k_\nu$ for $\nu = 1, \ldots, r$, then the intersection
$\bigcap\limits^r_{\nu = 1} A_\nu$ will "usually" be a flat of
dimension $d-\sum^r_{\nu=1} k_\nu$ (when $\sum_{\nu=1}^r k_\nu \le
d)$, or $\emptyset $ (when $\sum_{\nu=1}^r k_\nu > d).$ The dimension
of the  empty set $\emptyset$ is, by definition, $-1$.  We {\bf
always} have:

{\bf either} $\dim\bigcap^r_{\nu=1} A_\nu \ge d -
\sum^r_{\nu=1} k_\nu$,\:\:{\bf or}\: $\bigcap^r_{\nu=1} A_\nu = \emptyset$.

In view of these observations we define:
\begin{definition}  A finite set  $S\subset R^d$ is in {\bf strong
general position (SGP)} if:

(a)  $S$ is in general position, i.e., every subset of $S$ of size
$\leq d+1$ is affinely independent or, equivalently, $\dim \aff
F=\min(d, |F| - 1)$  for all subsets $F\subseteq S$.

(b) For any collection $\{ F_1, \ldots, F_r\}$ of $r$ pairwise
disjoint subsets of $S \quad (1 \le r \le |S|)$:
\begin{equation}
\label{equ:Det_SGP}
 d-\dim\bigcap^r_{\nu=1} \aff F_\nu = \min(d+1,
\sum^r_{\nu=1}(d-\dim\aff F_\nu)).
\end{equation}
\end{definition}

\begin{figure}[h]%
\includegraphics[width=0.5\columnwidth]{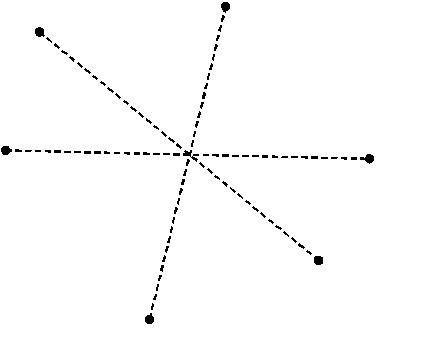}%
\caption{These six points are in general position but not in SGP}%
\end{figure}

Remark: Condition (\emph{a}) in the definition above follows from
condition (\emph{b}).  In fact, if $S$ is not in general position, i.e.,
if $S$ has an affinely dependent subset of size $\le d+1$,
consider a minimal affinely dependent subset $F$ of $S$.  Assume
$|F| = k, \; 3 \le k \le d+1$.  Then $\dim\aff F=k-2$.  The set
$F$ admits a Radon partition $F=A\cup B, \; A \cap B= \emptyset, \;
\conv A \cap \conv B \neq \emptyset$.  Assume $|A| = a (< k), |B| =
k - a(< k)$.  The sets $A, B$ are affinely independent.  Thus
$\aff A\cap \aff B \neq \emptyset$, even though
\begin{align*}
&d-\dim \aff A +d-\dim \aff B \\
&= d - (a-1) + d-(k-a-1)\\
&=2d-k+2\ge 2d-(d+1) + 2 = d+1.
\end{align*}

Our next aim is to show that if $S$ is a finite subset of $\bbr^d$
in general position, then we have to check only a small fraction
of the conditions listed in (\emph{b}) above in order to determine whether
$S$ is in SGP.  We shall do this in five steps.
The final result is stated as Theorem $\ref{thm:SGP_result}$ below.

(A) Suppose $S\subset \bbr^d$ is finite and in general position.  Then $S$ is in
SGP iff $(\ref{equ:Det_SGP})$ holds for any collection $F_1,\ldots, F_r$ of pairwise disjoint subsets of $S$ that satisfy $1\le |F_\nu | \le d$
for $\nu=1, 2,\ldots, r$.

\begin{proof}  If $F_\nu = \emptyset $ for some $1 \le \nu \le r$,
then $(\ref{equ:Det_SGP})$ holds automatically.  (Both sides of the equality are
$d+1$.)

If $|F_{\mu}| \ge d+1$ for some $1\le \mu \le r$, then $\aff F_{\mu} =\bbr^d$,
and removing $F_{\mu}$ from the collection does not affect the
intersection $\bigcap^r_{\nu = 1} \aff F_\nu $ on the left-hand side,
nor the sum $\sum^r_{\nu = 1} (d-\dim\aff F_\nu)$ on the right-hand
side.  Thus $(\ref{equ:Det_SGP})$ holds for the given collection $\{F_1, \ldots,
F_r\}$ iff it holds for the subcollection obtained by removing all
$F_\nu$'s with $|F_\nu| > d$.  (If all $|F_\nu|$'s are $> d$, then both
sides of $(\ref{equ:Det_SGP})$ are $0$.)
\end{proof}

(B) If $S\subset \bbr^d$ is finite and in general position, and $F_1,\ldots, F_r$ are
pairwise disjoint subsets of $S$ that satisfy $1\le |F_\nu|\le d$
for all $\nu$ and $\sum^r_{\nu=1} |F_\nu| = m$, then the equality
$(\ref{equ:Det_SGP})$ is equivalent to the condition:

\begin{equation}
\label{equ:Condition_B_for_SGP}
\begin{cases}  \dim \bigcap^r_{\nu=1} \aff F_\nu =
m-T(d, r) &\text{if } m \ge T(d, r)\\
\qquad \bigcap^r_{\nu=1} \aff F_\nu = \emptyset &\text{if } m < T(d,
r).
\end{cases}\end{equation}
\begin{proof} $(\ref{equ:Det_SGP})$ is equivalent to:
\begin{align*} \dim\bigcap^r_{\nu=1} \aff F_\nu &= d-\min (d+1,\,
\sum^r_{\nu=1} (d-\dim\aff F_\nu))\\
&= d-\min (d+1,\, rd+r-\sum^r_{\nu=1} |F_\nu|)\\
&=\max (-1,\, m - r(d+1) +d)\\
&=\max(-1,\, m - T(d, r)).\end{align*}.
\end{proof}

(C) Suppose $S\subset \bbr^d$ is finite and in general position.  Then $S$ is in SGP
iff: for any $r$ pairwise disjoint subsets $F_1, \cdots, F_r$ of
$S$ with $1\le |F_\nu| \le d$ for all $\nu \; (2\le r \le |S|)$ and
$\sum^r_{\nu=1} |F_\nu|= m$:
\begin{equation}
\label{equ:Intrsection_when_m<=T(d,r)}
\bigg|\bigcap^r_{\nu=1}\aff F_\nu\bigg| =
\begin{cases} 1 &\text{if } m = T(d, r)\\
0 &\text{if } m < T(d, r).\end{cases}\end{equation} (Note that
$(\ref{equ:Intrsection_when_m<=T(d,r)})$ always holds for $r=1$: when $r=1,\; m = |F_1|, \;T(d, 1) = 1$ and $\dim\aff F_1 = |F_1| -1=m-T(d, 1).)$

\begin{proof}
The "only if" direction is clear: condition $(\ref{equ:Intrsection_when_m<=T(d,r)})$ is the
restriction of condition $(\ref{equ:Condition_B_for_SGP})$ to the case $m\le T(d, r)$.

For the "if" direction: assume $(\ref{equ:Condition_B_for_SGP})$ fails for some $m  >
T(d, r)$, i.e., there are some $r$ pairwise disjoint subsets $F_1,
\ldots, F_r$ of $S$, $1 \le |F_\nu| \le d$ for all $\nu$, such that $m
= \sum^r_{\nu = 1} |F_\nu|  > T(d, r)$, and $\dim \bigcap^r_{\nu=1}
\aff F_\nu\neq m - T(d, r)$.

Note that this can happen only for $r\ge 2$.

 Among all the "violations" $(F_1, \ldots, F_r)$
of $(\ref{equ:Condition_B_for_SGP})$ with $m (= \sum_{\nu=1}^r |F_\nu|)> T(d, r)$ (where $r$ is not fixed in advance),
choose one with $m$ as small as possible.  Then one
of the following holds:

Case I: $\bigcap^r_{\nu=1}\aff  F_\nu = \emptyset$.

Case II: $\dim \bigcap^r_{\nu=1} \aff F_\nu > m-T(d, r)$.

In case I, choose nonempty subsets $G_\nu \subseteq F_\nu$ such that
$\sum^r_{\nu=1} |G_\nu| = T(d, r)$. This is possible, since $T(d, r)
= (d+1) (r-1) + 1 \ge r$. The sets $G_1, \ldots, G_r$  violate
condition $(\ref{equ:Intrsection_when_m<=T(d,r)})$, since $\bigcap^r_{\nu=1}\aff G_\nu \subset
\bigcap^r_{\nu=1} \aff F_\nu = \emptyset$, even though
$\sum^r_{\nu=1}|G_\nu| = T(d, r)$.

In case II, choose an index $\mu$ such that $|F_\mu| > 1$.  Pick a
point $p \in F_\mu$ and define $G_\mu = F_\mu \smallsetminus \{p\}$, $G_\nu =
F_\nu$ for all $\nu \neq \mu$.  Then $\dim\aff G_\mu = \dim \aff F_\mu-1$,
hence $\aff G_\mu = H \cap \aff F_\mu$ for some hyperplane $H \subset
\bbr^d$. Therefore, $\bigcap^r_{\nu=1}\aff G_\nu = H\cap
\bigcap^r_{\nu=1}\aff F_\nu$. This implies:

either

\begin{equation}
\label{equ:Empty_Intersection}
\bigcap^r_{\nu=1}
 \aff G_\nu = \emptyset,
 \end{equation}
 or
\begin{equation}
\label{equ:Large_dim_intersection}
\dim\bigcap^r_{\nu=1}
 \aff G_\nu \ge -1+\dim \bigcap^r_{\nu=1}
 \aff F_\nu > m-1-T(d, r),
 \end{equation}
where $m-1=\sum^r_{\nu=1} |G_\nu|$.

 If $(\ref{equ:Empty_Intersection})$ holds, then we have a smaller violation of $(\ref{equ:Condition_B_for_SGP})$ if $m-1>T(d, r)$ (contrary to our choice of $(F_1, \ldots, F_r))$, or a violation of $(\ref{equ:Intrsection_when_m<=T(d,r)})$, if $m-1=T(d, r)$.

 If $(\ref{equ:Large_dim_intersection})$ holds, then again we have a smaller violation of $(\ref{equ:Condition_B_for_SGP})$ if $m-1 > T(d, r)$, or a violation of $(\ref{equ:Intrsection_when_m<=T(d,r)})$ if $m-1 = T(d, r)$.

\end{proof}

In the next step we discard conditions that relate to the case $m<T(d,r)$ and are not minimal.

If $F_1,\ldots , F_r$ are pairwise disjoint subsets of $S$, and
$|F_\nu| = d+1-\varepsilon_\nu$, for $\nu=1, \ldots, r$, then $m=\sum^r_{\nu=1} |F_\nu| =
r(d+1) -\sum^r_{\nu=1}\varepsilon_\nu$, whereas $T(d, r) = r(d+1) - d$.
Thus $m < T(d, r)$ iff $\sum^r_{\nu=1} \varepsilon_\nu> d$. If, for some proper
subset $R'$ of $\{ 1, \ldots, r\}$ of size $r'$, we have
$\sum_{\nu\in R'} \varepsilon_\nu > d$, then $m'= \sum_{\nu \in R'}|F_\nu| =
r' (d+1) - \sum_{\nu\in R'} \varepsilon_\nu < r'(d+1) - d = T(d, r')$. If $\bigcap_{\nu\in R'} \aff F_\nu = \emptyset$, then, a
fortiori, $\bigcap^r_{\nu=1}\aff F_\nu = \emptyset.$

This reduces the criterion for SGP to the following:

(D) Suppose $S\subset \bbr^d$ is finite and in general position, $|S| > d+1$.  Then
$S$ is in SGP iff:

For any $r$ pairwise disjoint subsets $F_1, \ldots, F_r$ of $S \;(2
\le r \le d+1)$, if $|F_\nu| = d+1-\varepsilon_\nu$,  $1\le \varepsilon_\nu \le d$ for
$\nu = 1, \ldots, r$, then $\bigcap^r_{\nu=1} \aff F_\nu$ is a singleton
if $\sum^r_{\nu=1}\varepsilon_\nu = d$, and empty if $d<\sum^r_{\nu=1} \varepsilon_\nu \le d + \min\{\varepsilon_\nu,\: 1
\le \nu \le r\}$.

Now comes the final reduction:

\begin{thm}
   \label{thm:SGP_result}
 Suppose $S\subset \bbr^d$ is in general position and $|S| > d+1$.  Then
$S$ is in SGP iff:

For any $r$ pairwise disjoint subsets $F_1, \ldots, F_r$ of $S \;(2
\le r \le d+1)$, if $|F_\nu| = d+1-\varepsilon_\nu$,  $1\le \varepsilon_\nu \le d$ for
$\nu = 1, \ldots, r$, then $\bigcap^r_{\nu=1} \aff F_\nu$ is a singleton
if $\sum^r_{\nu=1}\varepsilon_\nu = d$, and empty if 

either $\sum^r_{\nu=1} \varepsilon_\nu = d + 1$,

or $|S|=r(d+1) -\sum^r_{\nu=1} \varepsilon_\nu $, and  $d+1<\sum^r_{\nu=1} \varepsilon_\nu \le d + \min\{\varepsilon_\nu,\: 1
\le \nu \le r\}$.
\end{thm}

Remarks: (a)  $|S|=r(d+1) -\sum^r_{\nu=1} \varepsilon_\nu $ means just that $F_1 \cup \cdots \cup F_r = S$.

(b) If $|S| \geq d(d+1)$, then we can dispense with the second clause in Theorem $\ref{thm:SGP_result}$, and the condition becomes:

\begin{equation}
\bigg|\bigcap^r_{\nu=1}\aff F_\nu\bigg| =
  \begin{cases} 
       1 &\text{if }\sum^r_{\nu=1} \varepsilon_\nu = d \\
       0 &\text{if }\sum^r_{\nu=1} \varepsilon_\nu = d + 1.
  \end{cases}
\end{equation}

\begin{proof}
The "only if" direction is clear: the conditions in Theorems $\ref{thm:SGP_result}$ are just a subset of the conditions in (D) above.

As for the "if" direction:

Assume one of the conditions in (D) that is missing in Theorem $\ref{thm:SGP_result}$ is violated: $F_1, \ldots , F_r$  are $r$ pairwise disjoint subsets of $S$,  $(2 \le r \le d+1)$, $|F_\nu| = d+1-\varepsilon_\nu$,  $1\le \varepsilon_\nu \le d$ for $\nu = 1, \ldots, r$,  $d+1<\sum^r_{\nu=1} \varepsilon_\nu \le d + \min\{\varepsilon_\nu,\: 1 \le \nu \le r\}$,  $|S|>r(d+1) -\sum^r_{\nu=1} \varepsilon_\nu $ (i.e. $S \supsetneq   \cup^r_{\nu=1} F_\nu$ ) and still $\cap^r_{\nu=1}\aff F_\nu \neq \emptyset.$

Choose such a violation with $\sum^r_{\nu=1} \varepsilon_\nu $ as small as possible. Note that  $\sum^r_{\nu=1} \varepsilon_\nu \geq d+2$.  Choose an index $\mu,\,  1 \leq \mu \leq r$, with $ \varepsilon_\mu \geq 2$, and a point $q \in S \smallsetminus  \cup^r_{\nu=1} F_\nu $, replace $F_ \mu $ by $F'_\mu = F_\mu \cup \{ q \}$, and define $F' _\nu= F_\nu $ for all $\nu \neq \mu$. Now $|F'_\nu|=d+1-\varepsilon_\nu'$, where $\varepsilon'_\nu=\varepsilon_\nu$ for $\nu \neq \mu$, and $\varepsilon'_\mu=\varepsilon_\mu-1$. Clearly $ \cap^r_{\nu=1} \aff F'_\nu \supset \cap^r_{\nu=1} \aff F_\nu \neq \emptyset $. If $\sum^r_{\nu=1} \varepsilon'_\nu>d+1$ and $|S|>r(d+1)-\sum^r_{\nu=1} \varepsilon'_\nu$, then we have a violation of a condition in (D) that is missing in Theorem $\ref{thm:SGP_result}$ with $\sum^r_{\nu=1} \varepsilon'_\nu<\sum^r_{\nu=1} \varepsilon_\nu$, contrary to our earlier choice. If $\sum^r_{\nu=1} \varepsilon'_\nu=d+1$, or $|S|=r(d+1)-\sum^r_{\nu=1} \varepsilon'_\nu$, then we have a violation of one of the conditions in Theorem $\ref{thm:SGP_result}$.
\end{proof}

Conclusion: The following "recipe" states  explicitly what has to be checked in order to ascertain that a given list $(a_1, \ldots , a_n)$ of points in $\bbr^d$ consists of $n$ distinct points in SGP. In order to learn how to perform the various checks, the reader is advised to consult Section $\ref{sec:PointAreUsuallyInSGP}$ below.

Step I: Check that the given points $a_1, \ldots , a_n$ are (distinct and) in (ordinary) general position.

Step II: Consider collections $\mathcal{F}=\{F_1, \ldots,F_r \}$ of pairwise disjoint subsets of $\{ a_1, \ldots , a_n \}$. Assume $1 \leq |F_\nu| \leq d$ for all $1 \leq \nu \leq r$, say $|F_\nu|=d+1-\varepsilon_\nu$, where $1 \leq \varepsilon_\nu \leq d$. (To avoid duplication, you may assume that $\min \{ i: a_i \in F_\nu\}<\min \{ i: a_i \in F_{\nu+1} \}$ for $\nu = 1,2,\ldots,r-1$.) Denote by $m$ the total size $\sum_{\nu=1}^r |F_\nu|$ of $\mathcal{F}$ ($m=r(d+1)-\sum_{\nu=1}^r \varepsilon_\nu $).

(A) If $m=T(d,r)$ (i.e.,$ \sum_{\nu=1}^r \varepsilon_\nu=d$), check that $|\cap_{\nu=1}^r \aff F_\nu|=1$. This should be done for all $r$, $2 \leq r \leq \min \{  d, [\frac{n+d}{d+1}]\}$.

(B) If $m=T(d,r)-1$ (i.e.,$ \sum_{\nu=1}^r \varepsilon_\nu=d+1$), check that $\cap_{\nu=1}^r \aff F_\nu=\emptyset$. This should be done for all $r$, $3 \leq r \leq \min \{  d+1, [\frac{n+d+1}{d+1}]\}$.

(C) Define $ \varepsilon_0 =\min \{  \varepsilon_1, \ldots,  \varepsilon_d \}$. If $n=m \leq T(d,r)-2$ (i.e., $ \sum_{\nu=1}^r \varepsilon_\nu \geq d+2$), but  $ \sum_{\nu=1}^r \varepsilon_\nu-\varepsilon_0 \leq d$ (which implies $\varepsilon_0 \geq 2$), check that $\cap_{\nu=1}^r \aff F_\nu=\emptyset$. This should be done only for  $3 \leq r =\lceil \frac{n+d+2}{d+1} \rceil \leq  [\frac{d+2}{2}]$. In fact, clause (C) is applicable iff $d \geq 4$,  $3 \leq r \leq  [\frac{d+2}{2}]$ and $T(d,r)-[\frac{d}{r-1}] \leq n \leq T(d,r)-2$.

\section{Points are "Usually" in SGP}
\label{sec:PointAreUsuallyInSGP}

Let $X=(\vec{x_{1}},\ldots,\vec{x_{t}})$ be a sequence of  $t$ points in $\bbr^d$. Denote by $x_{k1},\ldots, x_{kd}$ the coordinates of  $\vec{x_k}$  ($k=1,2,\ldots,t$). We regard the $td$ quantities $x_{k \nu}$ ($1 \leq k \leq t, 1 \leq \nu \leq d  $) as real variables, and propose to find a non-zero polynomial $P=P_{t,d}$ in these variables, in such a way that the points  $\vec{x_{1}},\ldots,\vec{x_{t}}$ are (distinct and) in SGP, unless $P(\vec{x_{1}},\ldots,\vec{x_{t}})=0$. As we have seen in the preceding section, strong general position is a conjuction of a long list of conditions. For each condition (E) on the list we shall produce a non-zero polynomial $P_E$, such that the violation of condition (E) by the points $\vec{x_{1}},\ldots,\vec{x_{t}}$ will imply $P_E(\vec{x_{1}},\ldots,\vec{x_{t}})=0$. The polynomial $P_{t,d}$  promised above will be the product of all these polynomials $P_E$.

Denote by $M(X)$ the $(d+1) \times t$ matrix whose $k$-th column consists of the number $1$, followed by the coordinates of $\vec{x_k}$, i.e., ${1 \choose {\vec{x_k}}}=(1,x_{k1}, \ldots, x_{kd})^t$. For a subsequence $B$ of $X$ of length $b$ we denote by $M(B)$ the $(d+1) \times b$ submatrix of $M(X)$ that consists of the columns that correspond to points of $B$ only.

Let us start with the condition that the points of $X$ be distinct and in (ordinary) general position. If  $t=d+1$, this means that the points $\vec{x_{1}},\ldots,\vec{x_{d+1}}$ are affinely independent, i.e., that $\det M(X) \neq 0$, so the corresponding polynomial is just $\det M(X)$. If $t>d+1$, this means that each $d+1$ of the points $\vec{x_{1}},\ldots,\vec{x_{t}}$ are affinely independent, so the corresponding polynomial is the product of the determinants of all ${t \choose {d+1}}$ $ (d+1) \times (d+1)$ square submatrices of $M(X)$. If $t<d+1$,  then general position of the points of $X$ is the same as affine independence, so the  condition is: $\rank M(X)=t$. This means that $M(X)$ has at least one $t \times t$ square non-singular submatrix, and the corresponding polynomial is the sum of squares of the determinants of all $d+1 \choose t$ $t \times t$ square submatrices of $M(X)$.

When $t \leq d+1$, SGP is the same as (ordinary) general position, so we can stop here. Assume, from now on, that $t>d+1$. We assume that the points   $\vec{x_{1}},\ldots,\vec{x_{t}}$ are (distinct and) in general position, (otherwise, some polynomial we found already vanishes at  $\vec{x_{1}},\ldots,\vec{x_{t}}$), and proceed with the additional conditions, as they appear in (D) in Section $\ref{sec:SGP}$ above. (To be precise, we use the notation of (D) ($|F_\nu| = d+1-\varepsilon_\nu$,  $1\le \varepsilon_\nu \le d$ for $\nu = 1, \ldots, r$), but we do not use the reduction from (C) to (D), except for the fact that $2 \leq r \leq d+1$.)

Let $F_1, \ldots , F_r$ ( $2 \leq r \leq d+1$) be disjoint subsets of $\vec{x_{1}},\ldots,\vec{x_{t}}$,  $|F_\nu| =  d+1-\varepsilon_\nu$,  $1\le \varepsilon_\nu < d$ for $\nu = 1, \ldots, r$.

Case I: If $\sum_{\nu =1}^r  \varepsilon_\nu =d$, then $|\cap_{\nu=1}^r \aff F_\nu|=1$. Denote by $\vec{z}$ the unique point of  $\cap_{\nu=1}^r \aff F_\nu$. For each $\nu$, $1 \leq \nu \leq r$, $\vec{z}$ can be expressed as an affine combination (i.e., a linear combination with sum of coefficients 1) of the points of $F_\nu$. This expression is unique since $F_\nu$ is affinely independent. Thus ${1 \choose {\vec{z}}} = \sum \{ \lambda_i {1 \choose {\vec{x_i}}} :x_ i \in F_\nu \}$ for $\nu = 1,2, \ldots ,r$. We can eliminate the point $\vec{z}$ from this system of equations by writing 

\begin{equation}
\label{equ:equationsSystem}
\begin{cases}  \sum \{ \lambda_i {1 \choose {\vec{x_i}}} :x_ i \in F_\nu \}= \sum \{ \lambda_i {1 \choose {\vec{x_i}}} :x_ i \in F_{\nu+1} \}$  for $\nu = 1,2, \ldots ,r-1,\\
\sum \{ \lambda_i  :x_ i \in F_1 \}=1.
\end{cases}\end{equation}

Let us order the points $\vec{x_i}$ within each block $F_\nu$ by increasing order of the index $i$, and the union $\cup_{\nu=1}^r F_\nu$ by letting $F_\mu$ precede $F_\nu$ whenever $\mu < \nu$. (In order to avoid duplication we could index the blocks $F_\nu$ by increasing order of the smallest index of their elements, i.e., $\mu < \nu$ iff $\min \{ i : x_i \in F_\mu \} <\min \{ i : x_i \in F_\nu \} $.) Denote by $\Lambda$ the column of coefficients $\lambda_i$, ordered correspondingly. The equations $(\ref{equ:equationsSystem})$ can be written as: 

\begin{equation}
\label{equ:A_lambda}
A \cdot \Lambda = \left(  \begin {array}{c} 1\\0\\ \vdots \\0 \end{array} \right)
\end{equation}

Where $A$  is a square matrix of order $T(d,r)(=1+(d+1)(r-1))$ as illustrated below:

\begin{figure}[h]%
\includegraphics[width=0.65\columnwidth]{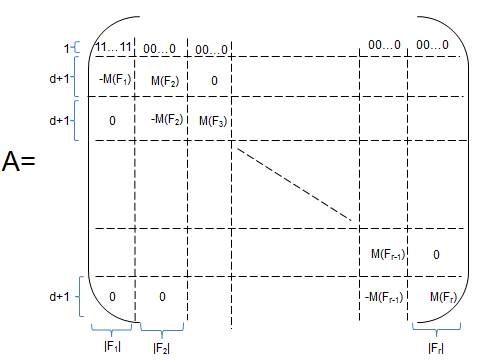}%
\caption{}%
\end{figure}

This system of (non-homogeneous) linear equations has a unique solution iff $\det A \neq 0$. Thus our polynomial is just $\det A$, regarded as a polynomial in the coordinates that appear as entries of $A$.

We still have to show that this polynomial is not identically $0$. Let $\vec{e_1} , \ldots ,\vec{e_d}$ be the standard orthonormal basis of $\bbr^d$. Recall that $|F_\nu|=d+1-\varepsilon_\nu$ ($\nu=1,2,\ldots,r$), where $1 \leq \varepsilon_\nu$ and $\sum_{\nu=1}^r \varepsilon_\nu=d$. For $\nu=1,2,\ldots,r$, define subspaces $W_\nu$ of $\bbr^d$ by: $W_\nu=\{ \vec{ x} \in \bbr^d : <\vec{e_i} ,\vec{ x}>=0$    for $ \sum_{\mu=1}^{\nu-1} \varepsilon_\mu<i\leq \sum_{\mu=1}^\nu \varepsilon_\mu \}$. Then $\dim W_\nu=d-\varepsilon_\nu=|F_\nu|-1$ and $\cap_{\nu=1}^r W_\nu=\{ \vec{0} \}$. For  $\nu=1,2,\ldots,r$, choose $F_\nu$ to be a set of $d+1-\varepsilon_\nu$ affinely independent points such that $\aff F_\nu=W_\nu$. Then  $\cap_{\nu=1}^r \aff F_\nu=\cap_{\nu=1}^r W_\nu=\{ \vec{0} \}$ is indeed a single point. Thus, for this choice of points $\det A \neq 0$, and therefore $\det A$, regarded as a polynomial in the coordinates of the vertices $\vec{x_i}$, is not identically $0$.

Case II: If $\sum_{\nu =1}^r  \varepsilon_\nu >d$ (i.e., if $\sum_{\nu=1}^r |F_\nu|<T(d,r)$) then $\cap_{\nu=1}^r \aff F_\nu = \emptyset$. This means that the system ($\ref{equ:A_lambda}$) $A \cdot \Lambda = \left(  \begin {array}{c} 1\\0\\ \vdots \\0 \end{array} \right)$ of linear equations, as described in case I, has no solution. Now $A$ is a rectangular $p \times q$ matrix, where $p=T(d,r)$, but  $q= \sum_{\nu=1}^r |F_\nu|<T(d,r)$. Denote by $A^+$ the augmented $p \times (q+1)$ matrix obtained by attaching the column $(1,0,\ldots,0)^t$ to $A$ ($q+1 \leq p$). 

Violation of the condition   $\cap_{\nu=1}^r \aff F_\nu = \emptyset$ means that the system  ($\ref{equ:A_lambda}$) $A \cdot \Lambda = (1,0,\ldots,0)^t$ does have a solution. This is equivalent to saying that the last column of $A^+$ is a linear combination of the first $q$ columns. This implies that $\rank A^+ \leq q$, which is equivalent to saying that all $(q+1) \times (q+1)$ submatrices of $A^+$ have zero determinant or (in view of the character  of the last column of $A^+$) that all $q \times q$ submatrices of $A$ that do not use the first row of $A$ have zero determinant.

 Denote by $A^-$ the rectangular matrix (of order $(d+1)(r-1)\times q$) obtained by deleting the first row of $A$, and let the polynomial $P$ be the sum of squares of all $q \times q$ subdeterminants of $A^-$. Violation of the condition   $\cap_{\nu=1}^r \aff F_\nu = \emptyset$ implies $P=0$.

We still have to show that this polynomial $P$ does not vanish identically. For any choice of vectors $\vec{x_1}, \ldots \vec{x_t}$ we have the quadruple equivalence: $P(\vec{x_1}, \ldots \vec{x_t})=0  \Longleftrightarrow \rank A^-<q  \Longleftrightarrow $ the columns of $A^-$ are linearly dependent $ \Longleftrightarrow $ the homogeneous system $A^- \cdot \Lambda =0$ of linear equations has a nonzero solution.

To complete the proof, we shall describe a particular substitution of vectors in $\bbr^d$ for the variable vectors  $\vec{x_1}, \ldots \vec{x_t}$, that will lead to a system $A^- \cdot \Lambda =\vec{0}$ whose only solution is $\Lambda =\vec{0}$ . Let  $\vec{u_0},\vec{u_1} \ldots \vec{u_d}$ be vectors in $\bbr^d$ whose only linear dependence (up to proportion) is $\sum_{i=0}^d \vec{u_i}=\vec{0}$. (Say, $\vec{u_i}=\vec{e_i}$ for $i=1,\ldots,d$ and $\vec{u_0}=-\sum_{i=1}^d \vec{u_i}$. Define $U=\{\vec{u_0},\vec{u_1} \ldots \vec{u_d} \}$. For $\nu=1,2,\ldots,r$ let $E_\nu$ be subsets of $U$ that satisfy: $|E_\nu|=\varepsilon _\nu$ (where $|F_\nu|=d+1-\varepsilon_\nu$), and $\cup_{\nu=1}^r E_\nu=U$. (Recall that $\sum_{\nu=1}^r \varepsilon_\nu \geq d+1$). Subtitute for the variable vectors $\vec{x_i} \in F_\nu$ (bijectively) the vectors $\vec{u_j} \in U \smallsetminus E_\nu$. (There is no need to substitute anything for variable vectors $\vec{x_i}$ that are not in $\cup_{\nu=1}^r F_\nu$, since they do not appear in $A$.) Note also that this substitution does not yield a set of $t$ points in SGP in $\bbr^d$: each point $\vec{u_j} \in U$ may appear up to $r-1$ times on the list  $\vec{x_1}, \ldots \vec{x_t}$). 

Now a solution $\Lambda$ of the homogeneous system of equations $A^- \cdot \Lambda =\vec{0}$ yields a point $\vec{z} \in \bbr^d$ that has $r$ representations 
\begin{equation}
\label{equ:representations_of_z}
\vec{z}= \sum_{\vec{u_j} \in U \smallsetminus E_\nu} \lambda_{j \nu}  \cdot \vec{u_j}\, \,  \text{ for }   \nu=1,\ldots,r  
 \end{equation}
where the sum of the coefficients is constant: $ \sum_{\vec{u_j} \in U \smallsetminus E_\nu} \lambda_{j\nu}= \sum_{\vec{u_j} \in U \smallsetminus E_{\nu+1}} \lambda_{j\,\,\nu+1}$ for $\nu = 1, \ldots , r-1$. The numbers $\lambda_{j\nu}$ ($ \nu=1,\ldots,r$, $\vec{u_j} \in U \smallsetminus E_\nu$) are the entries of the column vector $ \Lambda$. If $ \Lambda \neq \vec{0}$, then for some $\nu$, and some $\vec{u_{h}} \in U \smallsetminus E_\nu$,  $\lambda_{h\nu} \neq 0$. But $\vec{u_h} \in U=\cup_{\mu=1}^r E_\mu$, and therefore $\vec{u_h} \in E_\mu$ for some $\mu \neq \nu$. Consider the two representations: $\vec{z}=\sum \{  \lambda_{j\nu}\vec{u_j}:\vec{u_j} \in U \smallsetminus E_\nu\}=\sum \{  \lambda_{j\mu}\vec{u_j}:\vec{u_j} \in U \smallsetminus E_\mu\}$. They are different: $\vec{u_h}$ appears with a non-zero coefficient $\lambda_{h\nu}$ in the first one, but does not appear at all in the second one, since $\vec{u_h} \in E_
\mu$. In both representations, the sum of coefficients is the same. But this is impossible: If $\vec{z}=\sum_{j=0}^d \zeta_j \vec{u_j}$, then the most general representation of $\vec{z}$ as a linear combination of $\vec{u_0},\vec{u_1}, \ldots, \vec{u_d}$ is $\vec{z}=\sum_{j=0}^d (\zeta_j+\alpha) \vec{u_j}$, $\alpha \in \bbr$, so different representations necessarily differ in the sum of coefficients.


\end{document}